\documentclass[twocolumn]{autart}    
\usepackage{amsmath,amssymb,amsfonts}
\usepackage{dsfont}
\usepackage{graphicx}
\usepackage{algorithm,algorithmic}
\usepackage{color}
\usepackage{cite}
\usepackage{subcaption}
\usepackage{epstopdf}
\usepackage{bm}
\usepackage{mleftright}
\usepackage{mathrsfs}
\usepackage{stmaryrd}
\usepackage{apacite}
\usepackage{booktabs}
\usepackage{graphicx}
\usepackage{array}
\usepackage{subfloat}

\newtheorem{assumption}{Assumption}
\newtheorem{theorem}{Theorem}

\newtheorem{lemma}{Lemma}
\newtheorem{corollary}{Corollary}

\newtheorem{remark}{Remark}

\graphicspath{{./fig/}}
\begin{document}

\begin{frontmatter}

\title{Distributed Prediction-Correction Algorithms for Time-Varying Nash Equilibrium Tracking \thanksref{footnoteinfo}}

\thanks[footnoteinfo]{This work was sponsored by the National Key Research and Development Program of China under No. 2022YFA1004700, by the National Natural Science Foundation of China under Grant No. 62103343 and 62003239, by Shanghai Sailing Program under Grant Nos. 20YF1453000, , and partially by Shanghai Municipal Science and Technology Major Project No. 2021SHZDZX0100.}

\author[First]{Ziqin Chen}\ead{cxq0915@tongji.edu.cn},    
\author[Second]{Ji Ma}\ead{maji@xmu.edu.cn},
\author[First,ost]{Peng Yi}\ead{yipeng@tongji.edu.cn}
and \author[First]{Yiguang Hong}\ead{yghong@iss.ac.cn}
\thanks[ost]{Corresponding author.}
\address[First]{School of Electronics and Information Engineering,
	Tongji University, Shanghai, P. R. China}
\address[Second]{School of Aerospace Engineering,
	Xiamen University, Xiamen, P. R. China}

\vspace{-1.5em}

\begin{abstract}
This paper focuses on a time-varying Nash equilibrium trajectory tracking problem, that is applicable to a wide range of non-cooperative game applications arising in dynamic environments. To solve this problem, we propose a distributed prediction correction algorithm, termed DPCA, in which each player predicts future strategies based on previous observations and then exploits predictions to effectively track the NE trajectory by using one or multiple distributed gradient descent steps across a network. We rigorously demonstrate that the tracking sequence produced by the proposed algorithm is able to track the time-varying NE with a bounded error. We also show that the tracking error can be arbitrarily close to zero when the sampling period is small enough. Furthermore, we achieve linear convergence for the time-invariant Nash equilibrium seeking problem as a special case of our results. Finally, a numerical simulation of a multi-robot surveillance scenario verify the tracking performance and the prediction necessary for the proposed algorithm. 
\end{abstract}
\end{frontmatter}
\vspace{-0.5em}
\section{Introduction}
\vspace{-0.5em}
Game theory has been extensively used in numerous applications, including social networks~\cite{social}, sensor networks~\cite{sensor} and smart grids~\cite{smart}. The Nash equilibrium (NE), an important satisfactory solution concept in noncooperative games, emerges to be of both theoretical significance and practical relevance~\cite{NE}. Thus, the NE seeking method is becoming a key instrument in solving noncooperative game problems. 

Recently, with the penetration of multi-player systems, several distributed approaches to achieve the NE seeking have been proposed~\cite{staticNE1,staticNE2,staticNE3,staticNE4}. Such approaches, however, require time-invariant cost functions and NEs. In fact, in a variety of dynamic real-world situations, such as real-time traffic networks, online auctions, and intermittent renewable generations, the cost functions and NEs are time-varying. As a result, both academia and industry begin to focus on the time-varying NE tracking problem, which is a natural extension of time-invariant NE-related tasks in dynamic environments~\cite{online3,su2021online}. 

A natural ideal for solving the NE tracking problem is to treat it as a sequence of static problems and solve them one by one using  time-invariant NE seeking methods. However, these time-invariant methods necessitate iterative multiples or even infinite rounds for each static problem, and thus are not suitable for real-time tracking situations. To fix this drawback, a solvable approach is to predict the next NE at the current instant, and use the predicted variables as the initial value when solving the static problem at the next instant. Intuitively, only if the predicted model is properly constructed to bring the predicted variables close to the next NE, it is possible to achieve satisfactory tracking performance. This idea is called the prediction-correction technique, and it has been used in solving centralized time-varying optimization problems~\cite{simonetto2016}. The concept of changing initial values for iterative algorithms also motivates ADMM techniques in distributed time-varying optimizations~\cite{ADMM1}, which are to approximately solve each static optimization problem and use its output as the initial value for solving the next one. So far, neither prediction-correction techniques nor ADMM techniques have been used in solving distributed NE tracking problems to our knowledge. Another approach to solving NE tracking problems is found in the continuous-time regime~\cite{yemaojiao,multicluster}, which relies on ODE solutions and may not be implemented using digital computation.  In this paper, we focus on discrete-time algorithms.

The objective of this paper is to solve a time-varying NE tracking problem in a distributed manner with a prediction correction algorithm (DPCA). Specifically, a prediction step is proposed to predict the evolution of the NE trajectory using historical information, followed by a correction step to correct predictions by solving a time-invariant NE seeking problem obtained at each time. The contributions are summarized as follows.

1) To the best of our knowledge, our work is the first to study a distributed discrete-time algorithm for a time-varying NE tracking with provable convergence. 

2) We prove that the proposed DPCA can track time-varying NEs with a bounded error that is in inverse proportion to the number of correction steps. This result provides a quantitative analysis of the trade-off between tracking accuracy and algorithm iterations. Moreover, our proof demonstrates that the minimum number of correction steps can be one, thereby meeting the requirement of real-time tracking. 

3) Furthermore, we prove that with a small enough sampling period, the tracking error can be arbitrarily close to zero. For a time-invariant NE seeking problem, a special case of our formulation, we show that the tracking error linearly converges to zero.

4) A multi-robot surveillance scenario is used to test the performance of our algorithm with only one correction step. In addition, we compare the DPCA to no-predictor algorithms in terms of tracking error. The result verifies the significance of the introduced prediction step in the NE tracking and convergence.

This paper is organized as follows. The distributed NE trajectory tracking problem is presented in Section 2. The detailed algorithm design is proposed in Section 3 and the main results are established in Section 4. The result verification is investigated in Section 5. Finally, concluding remarks are given in Section 6.

\textit{Notation:}
Let $\mathbb{R}^{n}$ be the $n$-dimensional Euclidean space. For a scalar $a\in{\mathbb{R}}$, $\lceil a \rceil$ is the smallest integer not smaller than $a$. For $x\in{\mathbb{R}^{n}}$, denote its European norm by $\|x\|$ and its transpose by $x^{T}$. $[x_{i}]_{i\in{\mathcal{I}}}$ stacks the vector $x_{i}$ as a new column vector in the order  of the index set $\mathcal{I}$. For matrices $A$, the Kronecker product is producted by $\otimes$ and the smallest eigenvalue is denoted by $\lambda_{\min}(A)$. Denote by $0_{n},~1_{n}\in{\mathbb{R}^{n}}$, and $I_{n}\in{\mathbb{R}^{n\times n}}$ the vectors of all zeros and ones, and the identical matrix, respectively. For a differentiable function $J(x): \mathbb{R}^{n}\!\rightarrow\!\mathbb{R}$, denote $\nabla_{i}J(x)=\frac{\partial J(x)}{\partial x_i}$ with respect to $x_{i}$, where $x=[x_{i}]_{i\in \mathcal{I}}$.
\vspace{-0.5em}
\section{Problem Formulation}\label{sec:pre}
\vspace{-0.5em}
The time-varying game under consideration is denoted as $\Gamma^{t}=\{\mathcal{V},J^{t}_{i},x_{i}\}$, where $\mathcal{V}=\{1,\cdots,N\}$ denotes the set of $N$ players and $J_{i}^{t}: \mathbb{R}^{N}\rightarrow\mathbb{R}$ is the private time-varying cost function of player $i$. In the game $\Gamma^{t}$, each player $i$ aims to determine his time-varying strategy $x_{i}(t)$ to track the optimal trajectory shown below.
\vspace{-0.5em}
\begin{equation}
	x_i^*(t)=\underset{x_i\in \mathbb{R}}{\arg \min } J_i^{t}\left(x_i, x_{-i}\right),~i\in{\mathcal{V}},\label{p1}
\end{equation}
where $x_{-i}\in{\mathbb{R}^{N-1}}$ represents the joint strategies of all players except player $i$. If for any $t\geq0$, there is 
\begin{flalign}
	J^{t}_{i}(x^*_{i}(t),x^*_{-i}(t))=\underset{x_{i}(t)\in{\mathbb{R}}}{\operatorname{inf}} J^{t}_{i}(x_{i}(t),x^*_{-i}(t)),\nonumber
\end{flalign}
the time-varying vector $x^*(t)=[x_{i}^*(t)]_{i\in{\mathcal{V}}}$ is known as the NE trajectory. The NE trajectory tracking problem is a common issue in dynamic environments~\cite{multicluster}. Take a practical example, consider a target protecting problem in $\mathbb{R}_{+}^2$ by $N$ players in order to block $N$ intruders with locations $p_{i}^{t}$ and protect a time-varying target $b_{t}$. To protect the target to the greatest extent possible, it is preferable to keep the weighted center $\sigma(x)\!=\!\sum_{i=1}^{N}x_{i}(t)/N,~x\!=\![x_{i}]_{i\in{\mathcal{V}}}$ of all players to track the target $b_{t}$. In this case, each player is constantly adjusting his position to minimize the time varying $J_{i}^{t}\!=\!\|x_{i}\!-p_{i}^{t}\|^2+\gamma_{1}\|x_{i}\!-b_{t}\|+\gamma_{2}\|\sigma(x)\!-b_{t}\|^2$, and the problem can be cast as a NE trajectory tracking problem~\eqref{p1}. 

To distributively solve the NE trajectory tracking problem, each player $i$ estimates all other players' strategies for determining his tracking strategy over a communication graph. Assume that the underlying communication graph $\mathcal{G}(\mathcal{V},A)$ is undirected and connected. The adjacent matrix is defined as $A=[a_{ij}]_{N\times N}$, if the communication link between $i$ and $j$ exists, $a_{ij}>0$; otherwise, $a_{ij}=0$.  The corresponding Laplacian matrix is $L=[l_{ij}]_{N\times N}$ with $l_{ij}=-a_{ij}$ if $i \neq j$, and $l_{ij}=\sum_{j\neq i}^{N}a_{ij}$ otherwise. The neighbors set of player $i$ is defined as $\mathcal{N}_{i}=\{j\in\mathcal{V}|a_{ij}>0\}$. For an undirected and connected graph $\mathcal{G}$, all eigenvalues of $L$ are real numbers and can be arranged by an ascending order $0=\lambda_{1}<\cdots\leq \lambda_{N}$. 

Some definitions about player $i$'s estimated strategy are provided below. $x^{i}=[x^{i}_{i};x^{i}_{-i}]\in{\mathbb{R}^{N}}$ represents that player $i$'s estimates of all players' strategies; $x_{i}^{i}\in{\mathbb{R}}$ represents player $i$'s estimate of his strategy, which is indeed his actual strategy, i.e., $x_{i}^{i}=x_{i}$; $x_{-i}^{i}\in{\mathbb{R}^{N-1}}$ represents player $i$'s estimates on all players' strategies but his own. By the estimated strategy vector $x^{i}$ and the sampling period $h=t_{k+1}-t_{k}$, each player $i$ treat the problem~\eqref{p1} as a sequence of static problems
\begin{flalign}
	&\min _{x^{i}\in\mathbb{R}^{N}}\quad  J_{i}^{t_{k}}(x^{i}),\nonumber\\
	&\quad \text{s.t}\quad~~~ x^{i}=x^{j},~k\in{\mathbb{N}},\label{problem2}
\end{flalign}
and then solve them one by one. However, the sampling period $h$ in the real-time tracking is unable to afford the time spent in exactly solving every static problems by using time-invariant NE seeking algorithms in~\cite{staticNE2,staticNE3,staticNE4}. Therefore, how to design a dynamic algorithm to ensure good tracking performance with few iterative rounds during $h$ is the problem studied in this paper. 

To proceed, we make some standard assumptions, which were also required in relevant game works ~\cite{yemaojiao,online3}.
\begin{assumption}\label{A1}
The first and second derivatives of the time-varying NE trajectory exist and are bounded, i.e., $\|\dot{x}^*(t)\|\leq c_{1}$ and $\|\ddot{x}^*(t)\|\leq c_{2}$.
\end{assumption}
\begin{assumption}\label{A2}
Define the preudogradient as $F^{t}(x)\!=\![\nabla_{i}J_{i}^{t}(x)]_{i\in{\mathcal{V}}}: \mathbb{R}^{N}\!\rightarrow\! \mathbb{R}^{N}$. The following statements are required for any $t\geq 0$. \\	
	i) $F^{t}(x)$ is $\mu$-strongly monotone with $\mu>0$, that is, for any $x,y\in{\mathbb{R}^{N}}$,
	$$\langle F^{t}(x)-F^{t}(y),x-y\rangle\geq \mu\|x-y\|^2.\nonumber$$
	ii) There exist some constants $\theta_{i} \geq 0$ such that for any fixed $x_{-i}\in{\mathbb{R}^{N-1}}$, we have that for all $x_{i},y_{i}\in{\mathbb{R}}$,
	$$\|\nabla_{i}J^{t}_{i}(x_{i},x_{-i})-\nabla_{i}J^{t}_{i}(y_{i},x_{-i})\|\leq \theta_{i}\|x_{i}-y_{i}\|.\nonumber$$
	iii) There exist some constants $\theta_{-i}\geq 0$ such that for any fixed $x_{i}\in{\mathbb{R}}$, we have that for all $x_{-i},y_{-i}\in{\mathbb{R}^{N-1}},$
	$$\|\nabla_{-i}J^{t}_{i}(x_{i},x_{-i})\!-\!\nabla_{-i}J^{t}_{i}(x_{i},y_{-i})\|\!\leq\! \theta_{-i}\|x_{-i}\!-\!y_{-i}\|.$$
\end{assumption}
\begin{remark}
Assumption~\ref{A1} ensures that the variation of the time-varying NE is bounded, which is critical for bounded tracking errors. Assumption~\ref{A2}-i) guarantees the uniqueness of the NE trajectory. Assumptions~\ref{A2}-ii) and iii) correspond to the smooth of $J_i^t(x)$  on $x\in\mathbb{R}^{N}.$  
\end{remark}
The augmented preudogradient is further defined as $\boldsymbol{F}^{t}(\boldsymbol{x})\!\!=\!\![\nabla_{i}J^{t}_{i}(x^{i})]_{i\in{\mathcal{V}}}: \mathbb{R}^{N^{2}}\!\!\rightarrow\!\! \mathbb{R}^{N}$ with $\boldsymbol{x}\!=\![x^{i}]_{i\in{\mathcal{V}}}\!\in\! \mathbb{R}^{N^{2}}$. Meanwhile, denote $\theta\!=\!\max_{i\in\mathcal{V}}\{\sqrt{(\theta_{i})^2+(\theta_{-i})^2}\}$. Under Assumption~\ref{A2}, the $\theta$-Lipschitz continuity of $\boldsymbol{F}^{t}(\boldsymbol{x})$ is given in Lemma 1, whose proof is similar to that of Lemma 2 in~\cite{staticNE4} and is thus omitted.
\begin{lemma}\label{L1}
Under Assumption~\ref{A2}, the augmented preudogradient $\boldsymbol{F}^{t}(\boldsymbol{x})$ is $\theta$-Lipschitz continuous in $\boldsymbol{x}\in \mathbb{R}^{N^{2}}$.
\end{lemma}
\vspace{-0.5em}
\begin{table}
	\renewcommand\arraystretch{1.2}
	\label{dsmb}
	\centering
	\footnotesize
	\begin{tabular}{p{0.95\linewidth}}
		\toprule[1.5pt]
		{\textbf{Algorithm 1} Distributed Prediction Correction.} \\
		\midrule[1pt]
		\textbf{Initialization:} Initialize any $x^{i}_{0}\in{\mathbb{R}^{N}}$. Determine the number of correction steps $\tau\in{\mathbb{N}^{+}}$, the constant $q>0$, the stepsize $\alpha>0$ and $r_{i}\in{\mathbb{R}^{N}}$, which is a vector with $i$th element being $\varepsilon>0$ and the rest being $0$. \\
		1: \textbf{For} times $k=0,1,2\cdots$ \textbf{do},\\
		2: Generate predicted variables using prior information
		\vspace{-0.5em}
		\begin{eqnarray}
			\vspace{-0.5em}
			\quad \quad x^{i}_{k+1|k}=\left\{
			\begin{aligned}
				&x^{i}_{k},~~~~~~~~~~~~~~~~~~~~~~~\text{if}~k<q,\\
				&\frac{q+1}{q} x^{i}_{k}-\frac{1}{q}x^{i}_{k-q},~~~\text{otherwise}.
			\end{aligned}
			\right.\label{yuce}
			\vspace{-0.5em}
		\end{eqnarray}
		3:  Update the local cost functions
		$J_{i}^{t_{k+1}}(x)$.\\
		4: Initialize corrected variables $\hat{x}^{i,0}_{k+1}=x^{i}_{k+1|k}.$\\
		5: \textbf{for} $s=0: \tau-1$ \textbf{do}\\
		\quad \quad Receive neighbor's strategies $\hat{x}^{j,s}_{k+1}$ and execute
		\vspace{-0.5em}
		\begin{eqnarray}
			\quad \quad \hat{x}^{i,s+1}_{k+1}&&=\hat{x}^{i,s}_{k+1}+\alpha\sum_{i=1}^{N}a_{ij}(\hat{x}^{j,s}_{k+1}-\hat{x}^{i,s}_{k+1})\nonumber\\
			&&\quad \quad \quad \quad  \quad \quad  \quad \quad -r_{i}\nabla_{i}J_{i}^{t_{k+1}}(\hat{x}^{i,s}_{k+1}),\label{correctstep}
		\end{eqnarray}
		\quad~\textbf{end for.} \\
		6. Configure tracking variables $x^{i}_{k+1}=\hat{x}^{i,\tau}_{k+1}$.\\
		7. \textbf{end for}\\
		\bottomrule[1.5pt]
	\end{tabular}
\end{table}
\section{Algorithm Design}
In this section, Algorithm 1 is designed for each player $i,~i\!\in\!{\mathcal{V}}$ to generate a sequence $x^{i}_{k},~k\!\in\!{\mathbb{N}}$ to track time-varying NEs $x^*(t_{k})$. Algorithm 1 involves two steps: prediction and correction, which are depicted below.

\textbf{The prediction step design:} We derive the prediction step~\eqref{yuce} by observing the dynamic of time-varying NE. Due to $F^{t}(x^*(t))=0_{N}$, the NE trajectory $x^*(t)$ satisfies the following nonlinear dynamical system
\vspace{-0.5em}
\begin{flalign}
	\dot{x}^*(t)=F^{t}(x^*(t)).\label{L11dynamic}
\end{flalign}
Then we employ two methods to discretely estimate~\eqref{L11dynamic} and generate auxiliary prediction sequences $x_{k+1|k}\in{\mathbb{R}^{N}}$. Meanwhile, the prediction errors $\Delta_{k}=\|x_{k+1|k}-x^*(t_{k+1})\|$ of these two methods are compared.

i) The first method is to set the predicted variable as the NE $x^*(t_{k})$ at the last time instant $t_{k}$, in other words, 
\vspace{-0.5em}
\begin{flalign}
	x_{k+1|k}=x^*(t_{k}).
\end{flalign}
In this case, the prediction error is measured as
\vspace{-0.5em}
\begin{flalign}
	\Delta_{k}&=\|x^*(t_{k})-x^*(t_{k+1})\|,\nonumber\\
	&=\|\int_{t_{k}}^{t_{k+1}}\dot{x}^*(t)d\tau\| \leq c_{1}h.\label{D1}
\end{flalign}
ii) The second method involves a $q$-step interpolation method. In other words, 
\begin{flalign}
	\bar{x}_{k+1|k}=\frac{q+1}{q}x^*(t_{k})-\frac{1}{q}x^*(t_{k-q}),~k\geq q. \label{lemma12}
\end{flalign}
Recalling the dynamic~\eqref{L11dynamic}, the prediction error $\bar{\Delta}_{k}=\|\bar{x}_{k+1|k}-x^*(t_{k+1})\|$ is then computed as follows,
\begin{eqnarray}
	\bar{\Delta}_{k}&&=\left\|\frac{1}{q}x^*(t_{k})-\frac{1}{q}x^*(t_{k-q})-\int_{t_{k}}^{t_{k+1}}F^{\tau}(x^*)d\tau\right\|,\nonumber\\
	&&\leq \frac{1}{q}\sum_{r=0}^{q-1}\int_{t_{k}}^{t_{k+1}}\left\|F^{\tau}(x^*)-F^{\tau-(q-r)h}(x^*)\right\|d\tau. \label{yuce1}
\end{eqnarray}
Using the Mean Value Theorem and Assumption~\ref{A1} yields
\begin{eqnarray}
	\|F^{\tau}(x^*)-F^{\tau-(q-r)h}(x^*)\|&&=\left\|\frac{dF^{t}(x^*)}{dt}\Big{|}_{t=\xi}(q-r)h\right\|,\nonumber
\end{eqnarray}
where $\xi\in(\tau-(q-r)h,\tau)$. Assumption 1 ensures that
\begin{flalign}
\bar{\Delta}_{k}\leqslant \frac{c_{2}(q+1)h^{2}}{2}.\label{errorbound}
\end{flalign}
If the sampling period $h$ is small enough, then $\bar{\Delta}_{k}\leq {\Delta}_{k}$, which indicates that the second method outperforms the first. In fact, the sampling period $h\ll 1$ in the dynamic environment. Hence, it motivates the design of prediction in~\eqref{yuce}, while ensuring that predicted variables $x^{i}_{k+1|k}$ are close to $x^*(t_{k+1})$ at the time instant $t_{k+1}$.

\begin{remark}
The $q$ step interpolation method ensure the flexibility of the prediction's design. When $q=1$, the interpolation method in~\eqref{yuce} is related to two steps.
\end{remark}

\textbf{The correction step design:} We employ a distributed gradient descent (DGD) algorithm~\eqref{correctstep} to correct errors produced in the prediction. In particular, the predicted variable $x^{i}_{k+1|k}$ that is close to $x^*(t_{k+1})$ is harnessed as the correction step's initial value. It differs from the no-predictor DGD algorithm, which uses any initial state to exact update. In this sense, Algorithm 1 may require fewer iterative rounds $\tau$ and take less computation time to achieve the same tracking performance as the no-predictor DGD algorithm.
\vspace{-0.5em}
\section{Main Result}
\vspace{-0.5em}
This section establishes the convergence of Algorithm 1. We first analyze the prediction error in Lemma~\ref{L2}, and then we examine the tracking error in Lemma~\ref{L3}. Following those, Theorem 1 rigorously characterizes the upper bound of the tracking error. Before we begin, we define stacked variables $\boldsymbol{x}_{k+1|k}\!=\![x_{k+1|k}^{i}]_{i\in{\mathcal{V}}}$, $\boldsymbol{x}_{k}\!=\![x^{i}_{k}]_{i\in{\mathcal{V}}}$ and $\boldsymbol{x}^*(t_{k})\!=\!1_{N}\otimes x^*(t_{k})$. The proofs of Lemmas 2 and 3 can be found in Appendix. 
\begin{lemma}\label{L2}
Suppose Assumptions~\ref{A1} and~\ref{A2} hold. For any $k\geq q$, the norm of prediction errors is upper bounded by
\begin{eqnarray}
	&& \|\boldsymbol{x}_{k+1|k}-\boldsymbol{x}^{*}(t_{k+1})\|\leq\frac{q+1}{q}\|\boldsymbol{x}_{k}-\boldsymbol{x}^{*}(t_{k})\|\nonumber\\
	&&\quad \quad \quad  +\frac{1}{q}\|\boldsymbol{x}_{k-q}-\boldsymbol{x}^{*}(t_{k-q})\|+\frac{c_{2}\sqrt{N}(q+1)}{2}h^2.\label{predictionerror}
\end{eqnarray}
\end{lemma}
\begin{remark}
Lemma 2 gives the prediction error when $k\geq q$, but not when $k<q$. The reason is that the prediction error analysis for the case $k<q$ is able to be included in the proof of Theorem 1.
\end{remark}

\begin{lemma}\label{L3}
Suppose Assumptions~\ref{A1} and~\ref{A2} hold. There exist $\varepsilon$ and $\alpha$ satisfying
\begin{equation}
\begin{cases}
\alpha\in(0,\frac{\lambda_{2}}{2\lambda_{N}^2}),\nonumber\\
\varepsilon\in\left(0,\min\left\{\frac{\alpha\lambda_{2}\mu N}{4\theta^2+2\theta\mu N +\mu^2},\frac{\mu}{2N\theta^2},\frac{N}{4\alpha \lambda_{2}\mu}\right\}\right),
\end{cases}\label{lemma311}
\end{equation}
such that the following matrix is positive
\begin{flalign}
	A_{\alpha,\varepsilon}=\left[\begin{matrix}
		\small
		\frac{\varepsilon\mu}{N}-\varepsilon^2\theta^2 & -\frac{\varepsilon\theta}{\sqrt{N}}\\
		-\frac{\varepsilon\theta}{\sqrt{N}} & \alpha \lambda_{2}-\alpha^2\lambda_{N}^2-\varepsilon^2\theta^2-\varepsilon\theta
	\end{matrix}\right],
\end{flalign}
and $\nu=2\lambda_{\min}(A_{\alpha,\varepsilon})\in(0,1)$. Then for any $k\in{\mathbb{N}}$, the norm of correction errors is upper bounded by 
\begin{flalign}
	\|\boldsymbol{x}_{k+1}-\boldsymbol{x}^{*}(t_{k+1})\| \leq \rho^{\tau}\|\boldsymbol{x}_{k+1|k}-\boldsymbol{x}^{*}(t_{k+1})\|,\label{correctionerror}
\end{flalign}
where the coefficient $\rho=(1-\nu)^{\frac{1}{2}}\in(0,1)$.
\end{lemma}

\begin{remark}
Lemma 3 implies that the correction step~\eqref{correctstep} provides a discrete time approach to time-invariant distributed NE seeking problems in~\cite{staticNE1,staticNE2,staticNE3,staticNE4}. Meanwhile, Lemma 3 proves the linear convergence rate $O(\rho^{\tau})$ of this discrete time approach.
\end{remark}
With Lemmas~\ref{L2} and~\ref{L3}, we obtain our main theorem.
\begin{theorem}
Suppose Assumptions~\ref{A1} and~\ref{A2} hold. For any positive integer $q>\lceil\frac{2\rho^{\tau}}{1-\rho^{\tau}}\rceil$, there are positive constants $\gamma\in(0,1)$ and $C_{1}$ satisfying
\begin{flalign}
&\frac{\rho^{\tau}(q+1)}{q}\frac{1}{\gamma}+\frac{\rho^{\tau}}{q}\frac{1}{\gamma^{q+1}}\leq 1,\label{C12}\\
&C_{1}\geq\|\boldsymbol{x}_{0}-\boldsymbol{x}^*(t_{0})\|+\frac{\sqrt{N}\rho^{\tau}hc_{1}}{1-\rho^{\tau}}.\label{C11}
\end{flalign}
Then Algorithm 1 ensures the norm of tracking errors to be upper bounded by
\begin{equation}
 \|\boldsymbol{x}_{k}-\boldsymbol{x}^*(t_{k})\|\leq C_{1}\gamma^{k}+C_{2}h^2,\label{theorem1}
\end{equation}
where the positive constant $C_{2}=\frac{c_{2}\sqrt{N}q(q+1)\rho^{\tau}}{2q-2\rho^{\tau}(q+2)}$.
\end{theorem}
\begin{pf}
We prove that there must be a $\gamma_{0}\in(0,1)$ such that~\eqref{C12} holds. Define a decaying function $g(\gamma): \mathbb{R} \rightarrow \mathbb{R}$ related to $\gamma\in(0,+\infty)$ as follows,
\begin{flalign}
	g(\gamma)=\frac{\rho^{\tau}(q+1)}{q}\frac{1}{\gamma}+\frac{\rho^{\tau}}{q}\frac{1}{\gamma^{q+1}}.\nonumber
\end{flalign}
By $q>\lceil\frac{2\rho^{\tau}}{1-\rho^{\tau}}\rceil$, we can easily obtain that
\begin{flalign}
g(1)=\frac{\rho^{\tau}(q+1)}{q}+\frac{\rho^{\tau}}{q}<1.\nonumber
\end{flalign} 
Thus, there is a constant $\gamma_{0}\in(0,1)$ satisfying~\eqref{C12}.

Next, we prove~\eqref{theorem1} by starting with the case $k< q$, and then complete the proof for any $k\in{\mathbb{N}}$. 

For any $k<q$, there is $\boldsymbol{x}_{k+1|k}=\boldsymbol{x}_{k}$. According to Lemma~\ref{L3} and the triangle inequality, we have
\begin{eqnarray}
&&\quad \|\boldsymbol{x}_{k+1}-\boldsymbol{x}^{*}(t_{k+1})\|\nonumber\\
&&\leq \rho^{\tau}\|\boldsymbol{x}_{k}-\boldsymbol{x}^{*}(t_{k})\|+\rho^{\tau}\|\boldsymbol{x}^{*}(t_{k})-\boldsymbol{x}^{*}(t_{k+1})\|,\nonumber\\
&&\leq (\rho^{\tau})^{k+1}\|\boldsymbol{x}_{0}-\boldsymbol{x}^{*}(t_{0})\|+\sum_{r=0}^{k}(\rho^{\tau})^{r+1}\nonumber\\
&&\quad \quad\quad \quad\quad \quad \quad \quad\quad \quad \|\boldsymbol{x}^{*}(t_{k+1-r})-\boldsymbol{x}^{*}(t_{k-r})\|.\label{T11}
\end{eqnarray}
The upper bound of the last term of \eqref{T11} is analyzed. Based on~\eqref{D1}, we know that $\|x^*(t_{k+1})-x^*(t_{k})\|\leq hc_{1}$. By the sequence summation, we rewrite~\eqref{T11} as
\begin{flalign}
	&\quad \|\boldsymbol{x}_{k+1}-\boldsymbol{x}^{*}(t_{k+1})\|\nonumber\\
	&\leq (\rho^{\tau})^{k+1}\|\boldsymbol{x}_{0}-\boldsymbol{x}^{*}(t_{0})\|+\frac{\sqrt{N}\rho^{\tau}(1-(\rho^{\tau})^{k+1})hc_{1}}{1-\rho^{\tau}}.\nonumber
\end{flalign}
By \eqref{C12} with $\rho^{\tau}<\gamma<1$, we have that for any $k<q$,
\begin{flalign}
	\|\boldsymbol{x}_{k}-\boldsymbol{x}^{*}(t_{k})\|&\leq \gamma^{k}\|\boldsymbol{x}_{0}-\boldsymbol{x}^{*}(t_{0})\|+\frac{\sqrt{N}\rho^{\tau}hc_{1}}{1-\rho^{\tau}}
	\leq C_{1}\gamma^{k}.\nonumber
\end{flalign}
At last, we prove~\eqref{theorem1} for any $k\in{\mathbb{N}}$ via the mathematical induction. Assume that $\|\boldsymbol{x}_{k'}-\boldsymbol{x}^*(t_{k'})\|\leq C_{1}\gamma^{k'}+C_{2}h^2$ for any $k'\leq k,~k\geq q$. We prove $\|\boldsymbol{x}_{k'+1}-\boldsymbol{x}^*(t_{k'+1})\|\leq C_{1}\gamma^{k'+1}+C_{2}h^2.$ By inputting~\eqref{predictionerror} into~\eqref{correctionerror}, we have
\begin{flalign}
&\quad \|\boldsymbol{x}_{k'+1}-\boldsymbol{x}^{*}(t_{k'+1})\|\nonumber\\
&\leq \frac{\rho^{\tau}(q+1)}{q}\|\boldsymbol{x}_{k}-\boldsymbol{x}^{*}(t_{k})\|\nonumber\\
&\quad  +\frac{\rho^{\tau}}{q}\|\boldsymbol{x}_{k-q}-\boldsymbol{x}^{*}(t_{k-q})\|+\frac{c_{2}\rho^{\tau}\sqrt{N}(q+1)}{2}h^2.\label{C13}
\end{flalign}
Integrating the initial argument $\|\boldsymbol{x}_{k'}-\boldsymbol{x}^*(t_{k'})\|\!\leq\! C_{1}\gamma^{k'}\!+\!C_{2}h^2$ for any $k'\leq k,~k\geq q$, the definition of $C_{2}$ and the inequality~\eqref{C12} into~\eqref{C13}, we obtain that for any $k\geq q$, $\|\boldsymbol{x}_{k'+1}-\boldsymbol{x}^{*}(t_{k'+1})\|\leq C_{1}\gamma^{k'+1}+C_{2}h^2.$
Combining $\|\boldsymbol{x}_{k}-\boldsymbol{x}^*(t_{k})\|\leq C_{1}\gamma^{k}$ for any $k\in[0,q)$, we conclude that for any $k\in{\mathbb{N}}$, the inequality \eqref{theorem1} holds.
\end{pf}
Theorem 1 shows that each sequence $\{x^{i}_{k}\}$ approaches a neighborhood of the NE trajectory $x^*(t_{k})$ as time passes. The exponentially convergence of Algorithm 1 is shown in the following corollary.
\begin{corollary}
Under the same conditions of Theorem 1, the sequence $\{\boldsymbol{x}_{k}\}$ generated by Algorithm 1 converges to $\boldsymbol{x}^*(t_{k})$ exponentially up to a bounded error as
	\begin{flalign}
\limsup_{k\rightarrow \infty}\|\boldsymbol{x}_{k}-\boldsymbol{x}^*(t_{k})\| \leq C_{2}h^2=O(h^2).\nonumber
	\end{flalign}
\end{corollary}
\begin{remark}
Corollary 1 describes the exponential convergence property of Algorithm 1. The convergence accuracy is determined by the values of $h$ and $C_{2}$ that is in inverse proportion to the correction steps $\tau\in{\mathbb{N}^{+}}$. Note that increasing the value of $\tau$ improves tracking accuracy while necessitating much more computation time. Thus, the number of correction steps $\tau$ can be adjusted to match the actual tracking bound. Furthermore, if the sampling period $h$ is small enough, then Algorithm 1 ensure that the sequence $\{\boldsymbol{x}_{k}\}$ is sufficiently close to the NE trajectory.
\end{remark}
\vspace{-0.5em}
\section{Example}
\vspace{-0.5em}
A robotic surveillance scenario~\cite{fangzhen} is used to validate the effectiveness of Algorithm 1, including a tracking error comparison between Algorithm 1 and the no-predictor algorithm. Consider five robots collaborating to protect a time-varying target from five intruders. Each intruder $i$ moves according to the following rule:
\begin{flalign}\label{pit}
	p_{i}^{t}=p_{i}^{c}+5\begin{bmatrix}
		\cos(t/20) \\
		\sin(t/10)
	\end{bmatrix},
\end{flalign}
where $p_{i}^{c}=[8+6(i-1);8+6(i-1)],~i\in\{1,\cdots,5\}$. The target moves via a similar law with~\eqref{pit}, so we set it to be $b_{t}=\sum_{i=1}^{5}p_{i}^{t}/10.$ The dynamic protection strategy of each robot is determined by the following cost function
\vspace{-0.5em}
\begin{flalign}
	f^{t}_{i}(x_{i},x_{-i})&=\frac{\gamma_{1}}{2}\|x_{i}(t)-p_{i}^{t}\|^2+\frac{\gamma_{2}}{2}\|x_{i}(t)-b_{t}\|^2\nonumber\\
	&\quad \quad \quad \quad \quad \quad +\frac{\gamma_{3}}{10}\|\sigma(x(t))-b_{t}\|^2,
\end{flalign}
where $x_{i}(t)\in{\mathbb{R}^{2}}$ denotes the robot's position, $\gamma_{1}=\gamma_{2}=0.1$ and $\gamma_{3}=0.875$. The aggregative variable $\sigma(x(t))=\sum_{i=1}^{5}x_{i}(t)/5$. The network graph is described by a ring. The parameters for Algorithm 1 are set to $q=400$, $\alpha=0.05$, $\varepsilon=0.01$, $h=0.1$ and $\tau=1$. The element of initial tracking sequences $\{x_{0}^{i}\}\in{\mathbb{R}^{10}}$ is generated at random within an interval $[5+5(i-1),10+5(i-1)]$.

Fig. 1 shows that for any $i\in\{1,\cdots,5\}$, every element of $x^{i}_{k}\in{\mathbb{R}^{10}}$ generated by Algorithm 1 can effectively reach consensus and track every element of the NE trajectory $x^*(t)$. Fig. 2 depicts that how the network constructs a smarter formation that protects the time-varying target $b_{t}$ better. As shown in Fig. 3(a), the tracking error $\|\boldsymbol{x}_{k}-\boldsymbol{x}^*(t_{k})\|$ converges to a constant term $7h^{2}=0.07$. Fig. 3 (b) demonstrates the importance of the prediction step in the NE trajectory tracking.
\begin{figure}
	\centering
	\subfloat[First dimension of $x_{i}^*(t)$.]{\includegraphics[height=3cm,width=4.2cm]{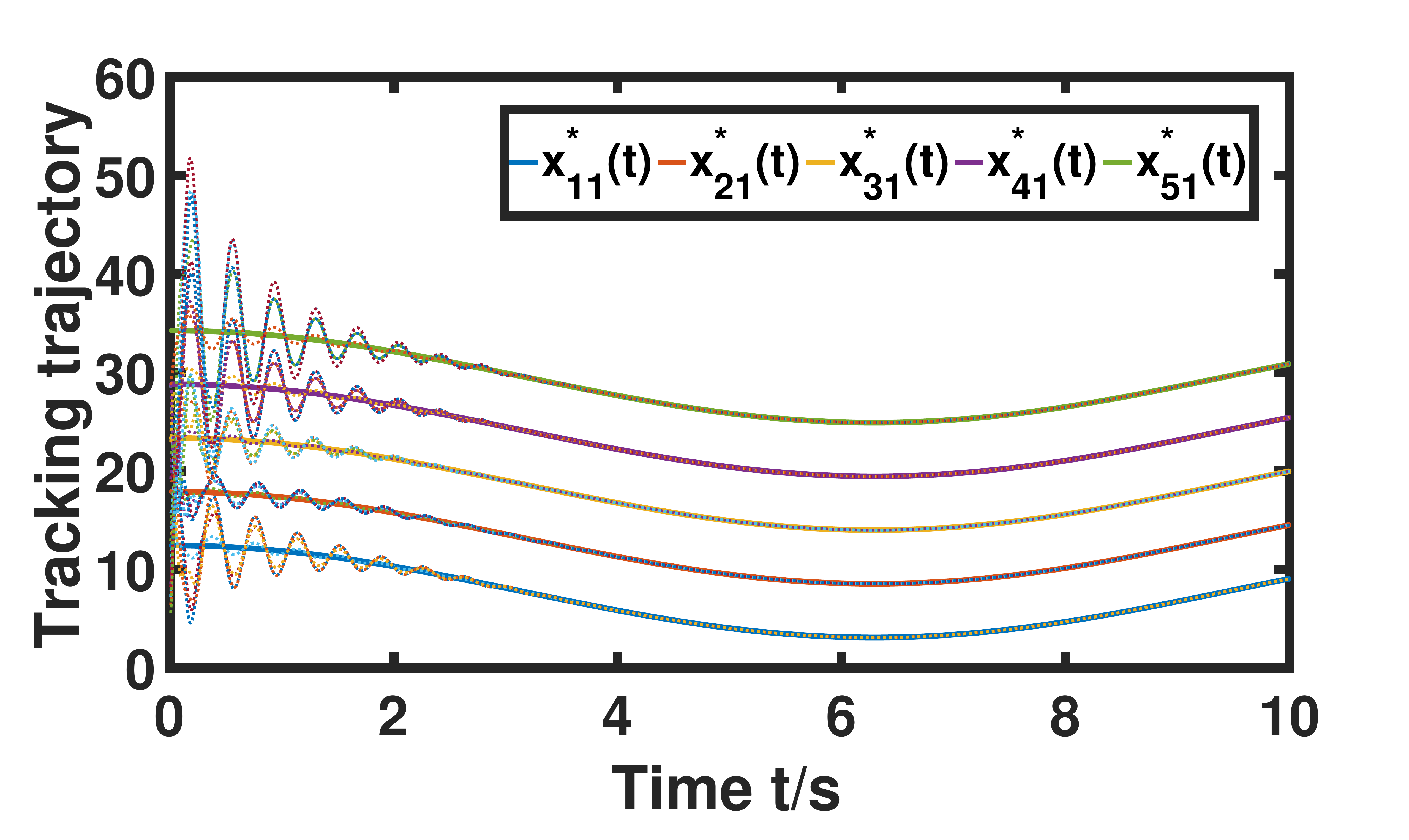}}
	\subfloat[Second dimension of $x_{i}^*(t)$.]{\includegraphics[height=3cm,width=4.2cm]{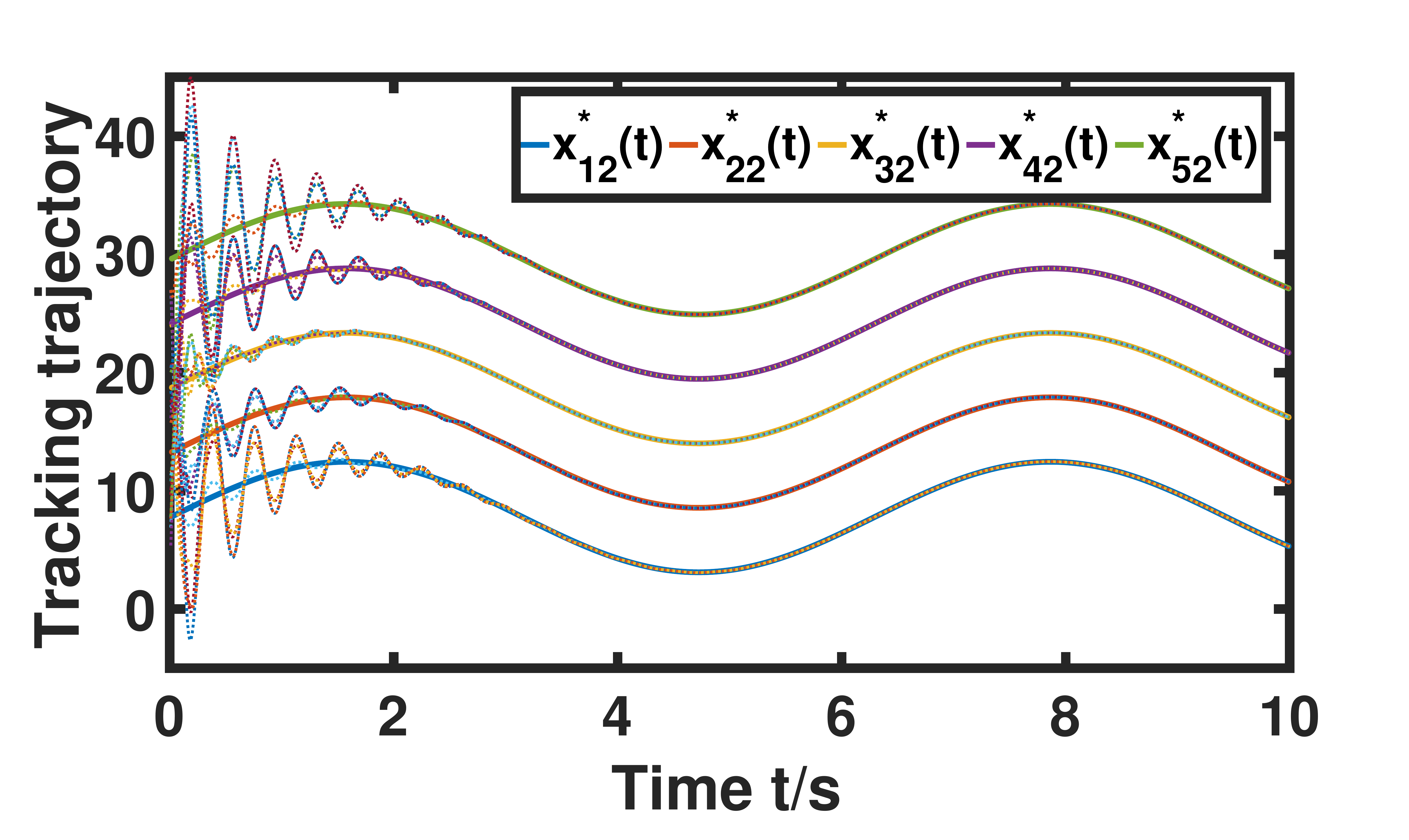}}
	\caption{Tracking performance of NE trajectory in terms of Algorithm 1. The solid lines represent $x^*(t)=[x_{i}^*(t)]_{i\in\{1,\cdots,5\}}$, where $x_{i}^*(t)=[x_{i1}^*(t);x_{i2}^{*}(t)]$. The dotted lines represent $x^{i}_{k}(t)\in{\mathbb{R}^{10}}$ in Algorithm 1.}
\end{figure}
\begin{figure}
	\centering
	\subfloat[t=0 s]{\includegraphics[height=2.5cm,width=4.2cm]{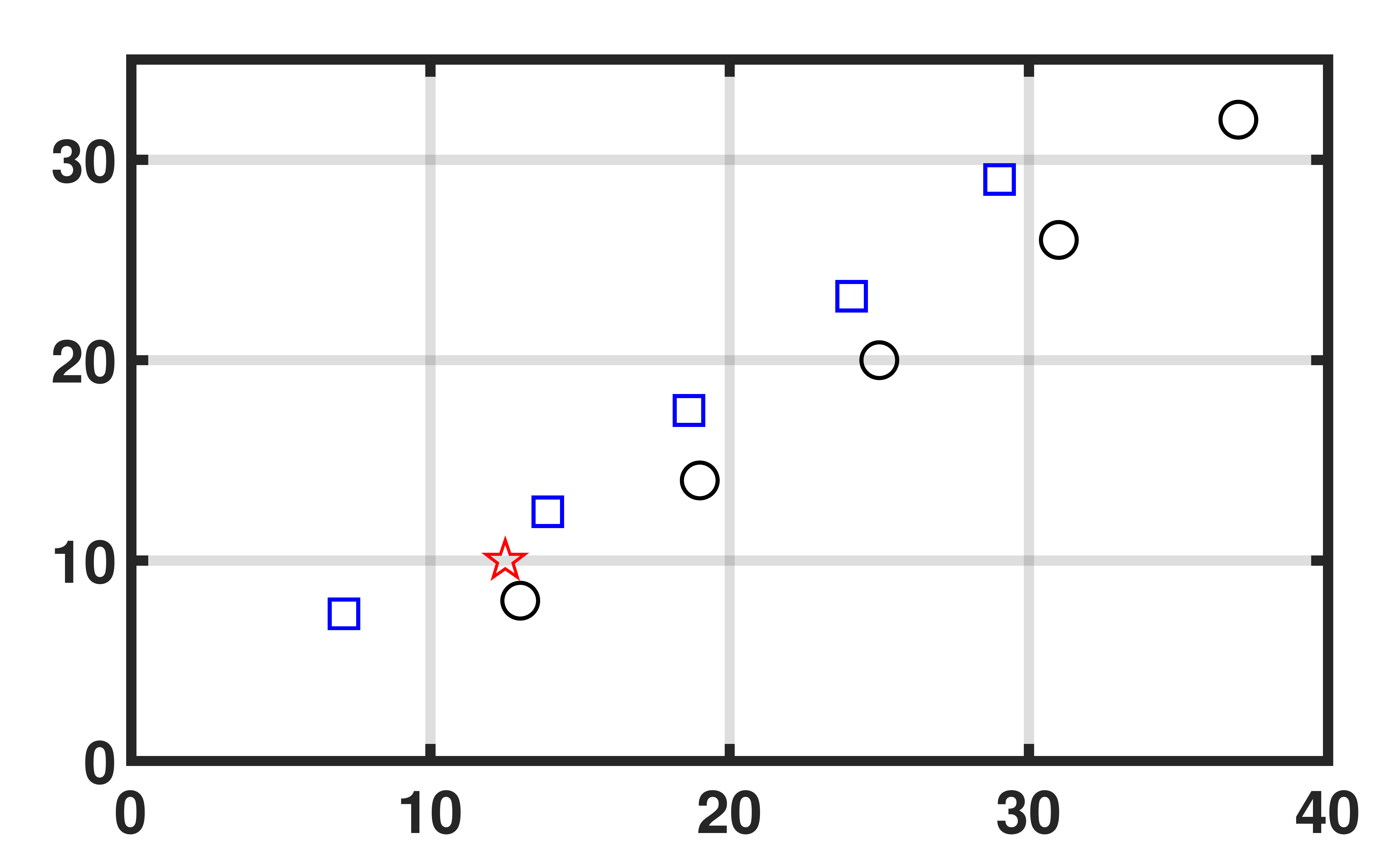}}
	\subfloat[t=3 s]{\includegraphics[height=2.5cm,width=4.2cm]{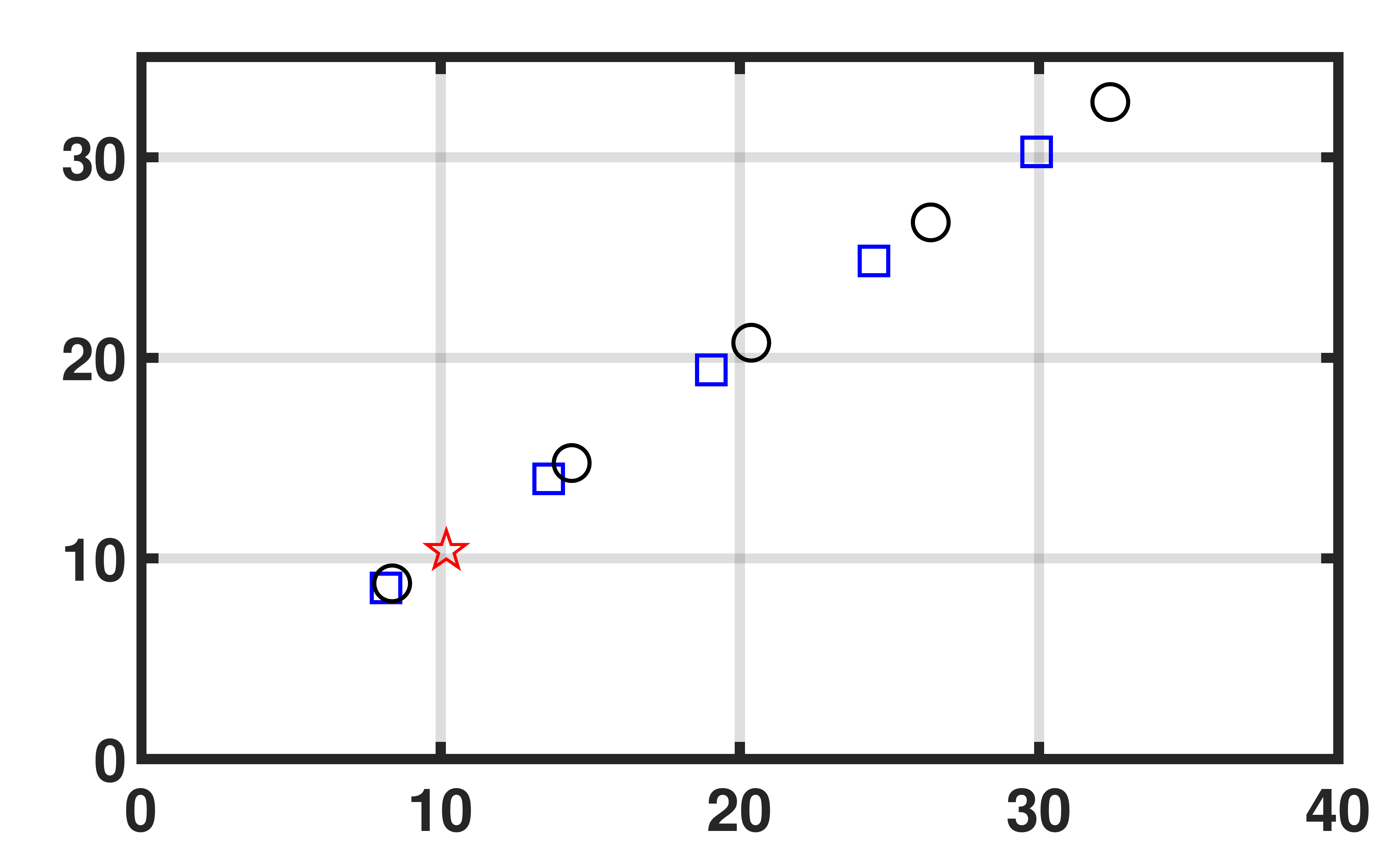}}	

	\subfloat[t=7 s]{\includegraphics[height=2.5cm,width=4.2cm]{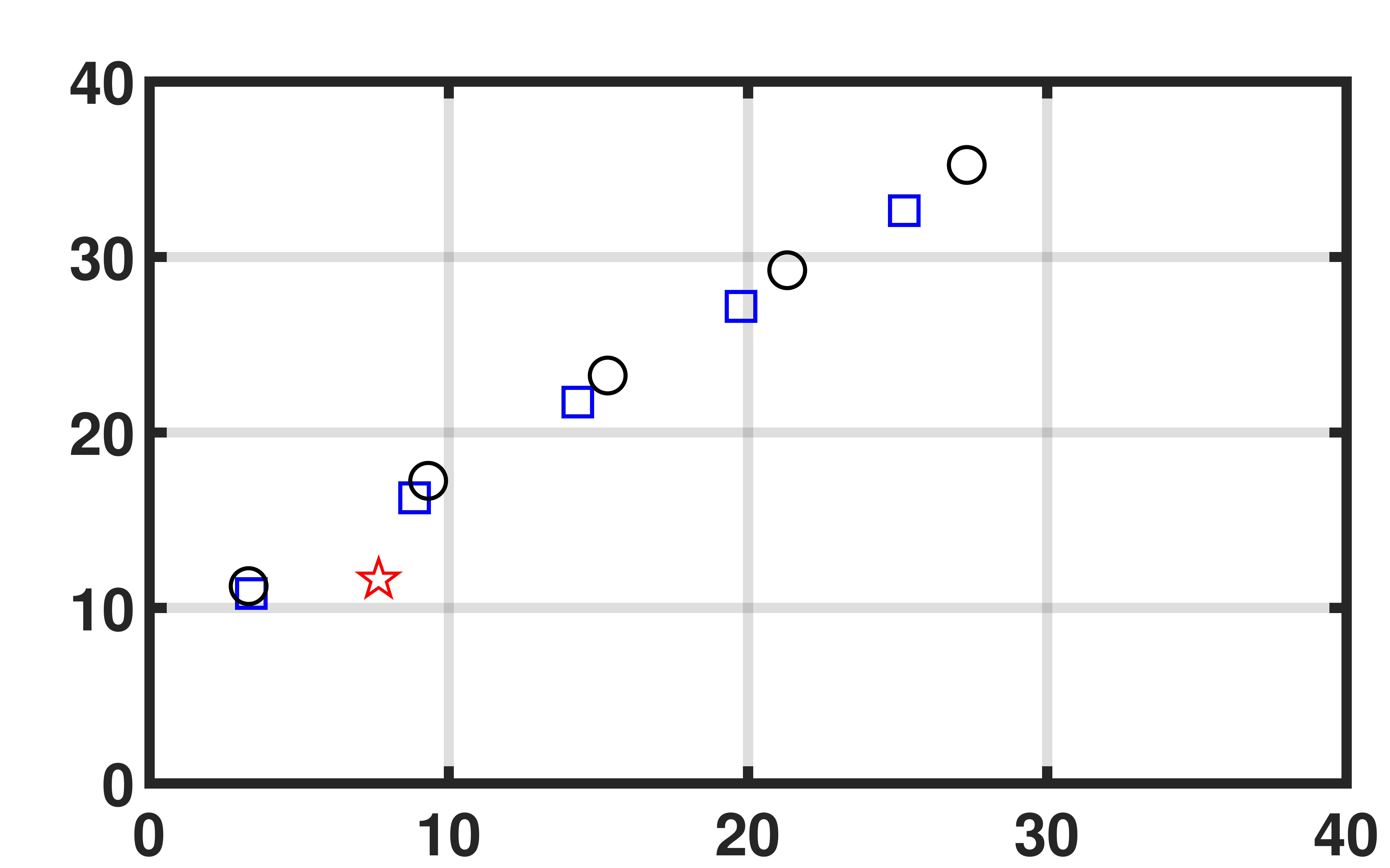}}
	\subfloat[t=10 s]{\includegraphics[height=2.5cm,width=4.2cm]{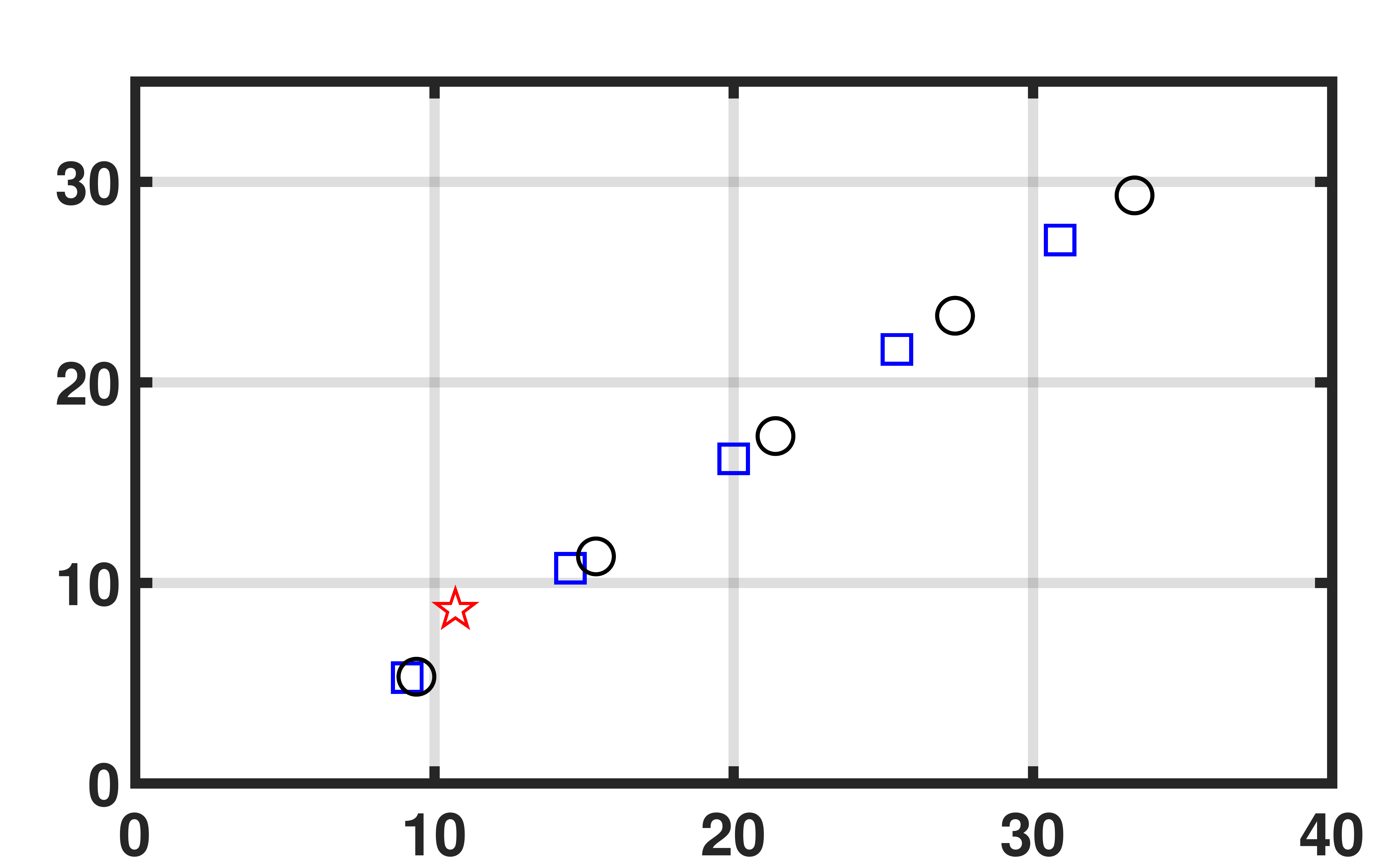}}
	\caption{The positions of robots, the time-varying intruders and targets at different time. The red star, circles and squares represent the target, robots and intruders, respectively.}
\end{figure}
\begin{figure}
	\centering
	\subfloat[DPCA's tracking error.]{\includegraphics[height=2.5cm,width=4.2cm]{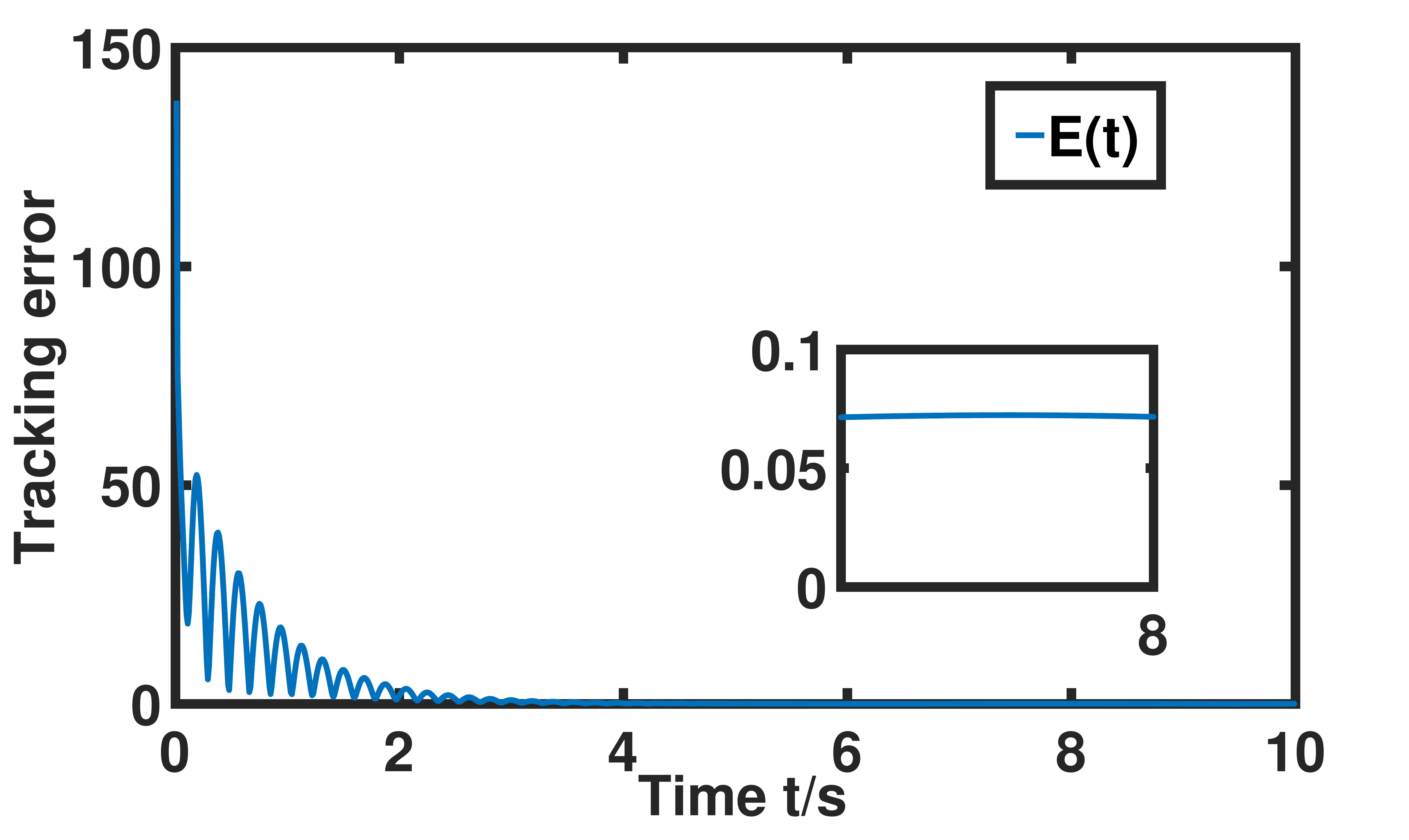}}
	\subfloat[The comparision between DPCA and no-predictor algorithm on tracking error.]{\includegraphics[height=2.5cm,width=4.2cm]{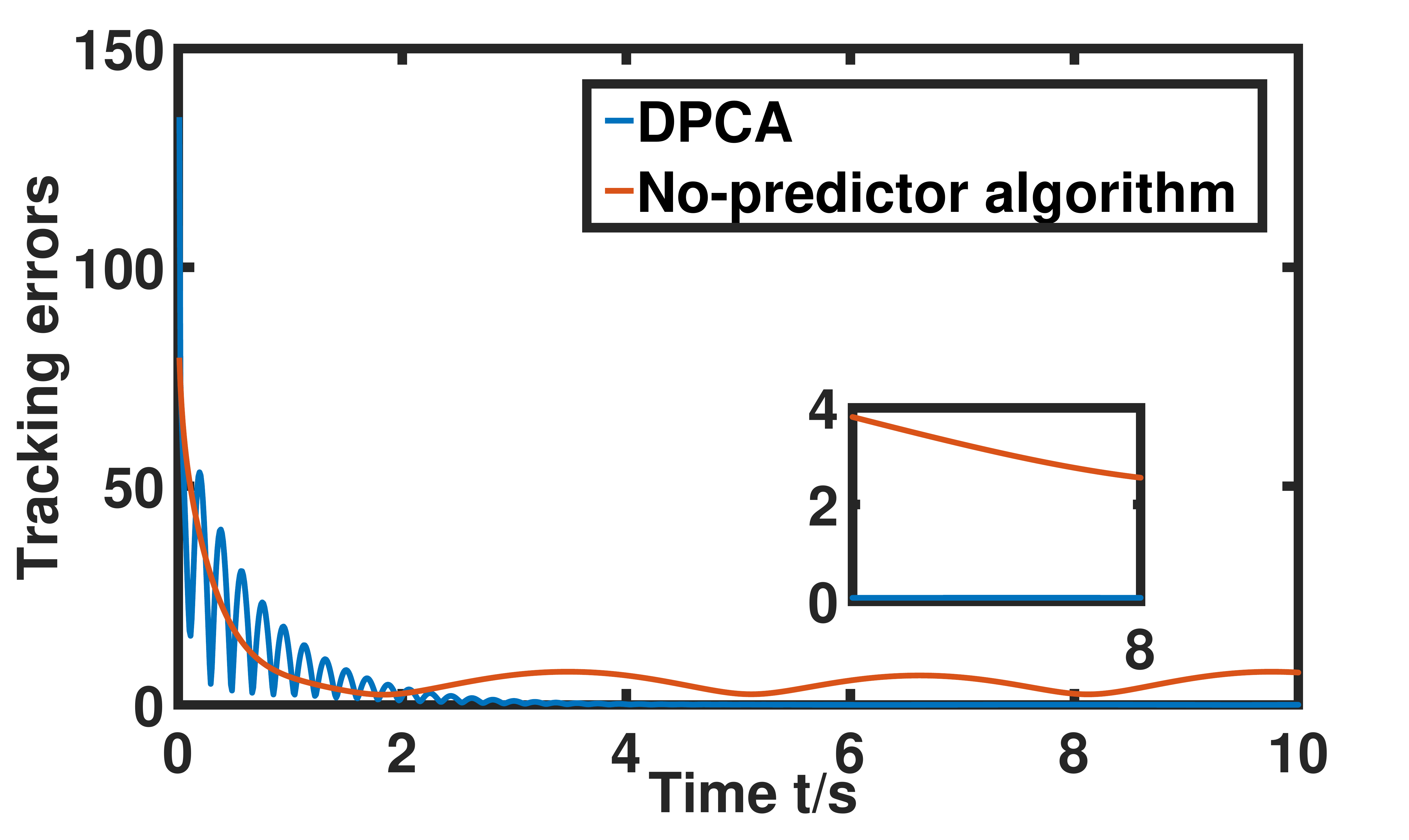}}
	\caption{The tracking error verification [cf. Theorem 1].}
\end{figure}
\section{Conclusion}
\vspace{-0.5em}
In this paper, we proposed a distributed prediction-correction algorithm (DPCA) to track the NE trajectory. By analyzing the coupling relationship between the prediction and the correction procedures, we rigorously proved that the tracking sequence generated by the DPCA can track the NE trajectory with a bounded error. This bound is proportional to the sampling period and inversely proportional to the number of correction steps. The results filled a gap in the study of NE tracking problems using discrete time algorithms. Moreover, the linear convergence rate for static games with time-invariant cost functions has been included in our results. The Future research of interests includes considering the communication efficient and privacy protected NE  tracking algorithm.  
\section{Appendix}
\vspace{-0.5em}
\subsection{Proof of Lemma 2}\label{plemma1}
In consideration of the error bound in~\eqref{errorbound}, we prove Lemma 2. By the prediction step in~\eqref{yuce}, we obtain
\begin{equation}
\begin{cases}
	&x^{i}_{k+1|k}=\frac{q+1}{q}x^{i}_{k}-\frac{1}{q}x^{i}_{k-q},~k\geq q,\\
	&x^*(t_{k+1})=\frac{q+1}{q}x^*(t_{k})-\frac{1}{q}x^*(t_{k-q})-\bar{\Delta}_{k}.\label{yuce2}
\end{cases}
\end{equation}
Then the prediction error can be computed as
\begin{flalign}
	&x^{i}_{k+1|k}-x^*(t_{k+1})\nonumber\\
	&\quad =\frac{q+1}{q}(x^{i}_{k}-x^*(t_{k}))-\frac{1}{q}(x^{i}_{k-q}-x^*(t_{k-q}))+\bar{\Delta}_{k}.\nonumber
\end{flalign}
Stacking the above equality, considering the norm of the resulting expression and applying the triangle inequality, Lemma 2 with~\eqref{predictionerror} holds.

\subsection{Proof of Lemma 3}
\vspace{-0.5em}
The proof Lemma 3 is divided into two parts: i) Prove that the matrix $A_{\alpha,\epsilon}$ is positive and the constant $\nu\in(0,1)$; ii) Prove the relationship between the tracking error and the prediction error in~\eqref{correctionerror}.

i) By $0<\varepsilon<\min\{\frac{\alpha\lambda_{2}\mu N}{4\theta^2+2\theta\mu N +\mu^2},\frac{\mu}{2N\theta^2},\frac{N}{4\alpha \lambda_{2}\mu}\}$, it yields
\begin{flalign}
\frac{\mu}{2N}>\varepsilon\theta^2~\text{and}~
\frac{\alpha\lambda_{2}}{2\varepsilon}> \frac{2\theta^2}{\mu N}+\theta+\frac{\mu}{2N} >\frac{2\theta^2}{\mu N}+\theta+\varepsilon\theta^2.\nonumber
\end{flalign}
Meanwhile, $\alpha\in(0,\frac{\lambda_{2}}{2\lambda^2_{N}})$ ensures $\frac{\alpha \lambda_{2}}{2\varepsilon}>\frac{\alpha^2\lambda_{N}^2}{\varepsilon}$, and thus,
\begin{flalign}
	\frac{\mu}{2N}\left[\frac{\alpha\lambda_{2}}{\varepsilon}-\frac{\alpha^2\lambda_{N}^2}{\varepsilon}-\varepsilon\theta^2-\theta\right] \geq \frac{\theta^2}{N^2}.\nonumber
\end{flalign}
By using again $\varepsilon\theta^2\leq\frac{\mu}{2N}$ and multiplying both side of the deriving results by $\varepsilon^2$, we have
\begin{flalign}
	(\frac{\varepsilon\mu}{N}-\varepsilon^2\theta^2)\left[\alpha\lambda_{2}-\alpha^2\lambda_{N}^2-\varepsilon^2\theta^2-\varepsilon\theta\right]-\frac{\varepsilon^2\theta^2}{N^2}>0,\nonumber
\end{flalign}
which is the determinant of $A_{\alpha,\varepsilon}$. Thus, $A_{\alpha,\varepsilon}>0$. By further considering $\varepsilon<\frac{N}{4\alpha \lambda_{2}\mu}$, it yields $\frac{\varepsilon\mu}{N}\alpha\lambda_{2}<\frac{1}{4}$ and 
\begin{flalign}
	(\frac{\varepsilon\mu}{N}-\varepsilon^2\theta^2)\left[\alpha\lambda_{2}-\alpha^2\lambda_{N}^2-\varepsilon^2\theta^2-\varepsilon\theta\right]-\frac{\varepsilon^2\theta^2}{N^2}<\frac{1}{4}.\nonumber
\end{flalign}
which indicates $\lambda_{\min}(A_{\alpha,\varepsilon})<\frac{1}{2}$. Thus, $\nu\in(0,1)$.

ii) Defining $\boldsymbol{L}=L\otimes I_{N}$ and $\hat{\boldsymbol{x}}^{s}_{k+1}=[\hat{x}^{i,s}_{k+1}]_{i\in \mathcal{V}}$, the correction step~\eqref{correctstep} can be written as a stacked form
\begin{flalign}
\hat{\boldsymbol{x}}^{s+1}_{k+1}-\hat{\boldsymbol{x}}^{s}_{k+1}=-\alpha \boldsymbol{L}\hat{\boldsymbol{x}}^{s}_{k+1}-R \boldsymbol{F}^{t_{k+1}}(\hat{\boldsymbol{x}}^{s}_{k+1}),\label{compact}
\end{flalign}
where $R=\text{diag}(r_{1},\cdots,r_{N})$. We first prove that at the equilibrium of system~\eqref{compact}, $\hat{x}^{i}_{k+1}$ reaches the NE of game $\Gamma^{t_{k+1}}$. If we set the left hand side of ~\eqref{compact} to zero, then
\begin{equation}
\alpha\boldsymbol{L} \hat{\boldsymbol{x}}^{*}_{k+1}+R \boldsymbol{F}^{t_{k+1}}(\hat{\boldsymbol{x}}^{*}_{k+1})=0_{N^2}.\label{kkt}
\end{equation}
Premultiplying both sides by $1^{T}_{N}\otimes I_{N}$ gives
\begin{flalign}
	&0_{N}=\alpha\left(1_N^{T} \otimes I_N\right)\left(L \otimes I_N\right) \hat{\boldsymbol{x}}^{*}_{k+1}\nonumber\\
	&\quad \quad \quad \quad \quad +\left(1_N^{T} \otimes I_N\right) R \boldsymbol{F}^{t_{k+1}}\left(\hat{\boldsymbol{x}}^{*}_{k+1}\right).\nonumber
\end{flalign}
Using $1_N^{T}L=0^{T}_{N}$ yields $0_{N}=\left(1_N^{T} \otimes I_N\right) R \boldsymbol{F}^{t_{k+1}}\left(\hat{\boldsymbol{x}}^*_{k+1}\right)$. By the definitions of $R$ and $\boldsymbol{F}^{t_{k+1}}$, $\boldsymbol{F}^{t_{k+1}}(\hat{\boldsymbol{x}}^*_{k+1})=0_{N}$, which further implies that $\boldsymbol{L}\hat{\boldsymbol{x}}^*_{k+1}=0_{N^2}$ from~\eqref{kkt}. There exist some $x^*(t_{k+1})\in{\mathbb{R}^{N}}$ such that $\hat{\boldsymbol{x}}^*_{k+1}=1_{N}\otimes x^*(t_{k+1})$ by the property of $L$ under the undirected connected graph $\mathcal{G}$. Thus, $\boldsymbol{F}^{t_{k+1}}(1_{N}\otimes x^*(t_{k+1}))=0_{N}$ and $F^{t_{k+1}}(x^*(t_{k+1}))=0_{N}$. In other words, $x^*(t_{k+1})$ is the NE of $\Gamma^{t_{k+1}}$ and $\hat{\boldsymbol{x}}^*_{k+1}=1_{N}\otimes x^*(t_{k+1})=\boldsymbol{x}^*(t_{k+1})$.

Next, we prove \eqref{correctionerror}. Before doing it, some coordinate transformations are required. Define $M=[M_{1};M_{2}]$ with $M_{1}=\frac{1_{N}}{\sqrt{N}}$ and $M_{2}\in{\mathbb{R}^{N\times{N-1}}}$, which satisfy $M_{2}^{T}M_{1}=0_{N-1}$, $M_{2}^{T}M_{2}=I_{N-1}$, and $M_{2}M_{2}^{T}=I_{N}-M_{1}M_{1}^{T}$. The unitary matrix $M$ can convert $L$ to a diagonal form $M^{T}LM=\text{diag}(0,\lambda_{2},\cdots,\lambda_{N})$. Denote   $\tilde{\boldsymbol{x}}^{s}_{k+1}=\hat{\boldsymbol{x}}^{s}_{k+1}-\boldsymbol{x}^*(t_{k+1})$ and perform the coordinate transformation as
\begin{eqnarray}
\bar{\boldsymbol{x}}^{s,1}_{k+1}=(M_{1}^{T}\otimes I_{N})\tilde{\boldsymbol{x}}^{s}_{k+1},~\bar{\boldsymbol{x}}^{s,2}_{k+1}=(M_{2}^{T}\otimes I_{N})\tilde{\boldsymbol{x}}^{s}_{k+1}.\label{255}
\end{eqnarray}
Then inputting~\eqref{compact} into~\eqref{255}, we obtain that
\begin{eqnarray}
\bar{\boldsymbol{x}}^{s+1,1}_{k+1}-\bar{\boldsymbol{x}}^{s,1}_{k+1}&&=-\left(M_1^{T} \otimes I_N\right) R\Theta^{s}_{k+1},\nonumber \\
\bar{\boldsymbol{x}}^{s+1,2}_{k+1}-\bar{\boldsymbol{x}}^{s,2}_{k+1}&&=-\alpha\left[\left(M_2^{T} L M_2\right) \otimes I_N\right] \bar{\boldsymbol{x}}^{s,2}_{k+1}\nonumber\\
&&\quad \quad\quad\quad\quad-\left(M_2^{T} \otimes I_N\right) R\Theta^{s}_{k+1},\label{311}
\end{eqnarray}
where $\Theta^{s}_{k+1}=\boldsymbol{F}^{t_{k+1}}\left(\hat{\boldsymbol{x}}^{s}_{k+1}\right)-\boldsymbol{F}^{t_{k+1}}\left(\boldsymbol{x}^*(t_{k+1})\right)$. 

Construct a Lyapunov function as follows,
\begin{flalign}
	V^{s}_{k+1}(\bar{\boldsymbol{x}}^{s,1}_{k+1},\bar{\boldsymbol{x}}^{s,2}_{k+1})=\|\bar{\boldsymbol{x}}^{s,1}_{k+1}\|^2+\|\bar{\boldsymbol{x}}^{s,2}_{k+1}\|^2.\nonumber
\end{flalign}
By the facts $(a+b)^{T}(a+ b)-b^{T}b=2a^{T}b+a^{T}a$ for any vectors $a,b$ and $\bar{\boldsymbol{x}}^{s+1,1}_{k+1}=\bar{\boldsymbol{x}}^{s+1,1}_{k+1}-\bar{\boldsymbol{x}}^{s,1}_{k+1}+\bar{\boldsymbol{x}}^{s,1}_{k+1}$, it yields
\begin{eqnarray}
	&&V^{s+1}_{k+1}-V^{s}_{k+1}=[\bar{\boldsymbol{x}}^{s+1,1}_{k+1}-\bar{\boldsymbol{x}}^{s,1}_{k+1}]^{T}[\bar{\boldsymbol{x}}^{s+1,1}_{k+1}-\bar{\boldsymbol{x}}^{s,1}_{k+1}]\nonumber\\
	&&
	\quad  +2[\bar{\boldsymbol{x}}^{s+1,1}_{k+1}-\bar{\boldsymbol{x}}^{s,1}_{k+1}]^{T}\bar{\boldsymbol{x}}^{s,1}_{k+1}+2[\bar{\boldsymbol{x}}^{s+1,2}_{k+1}-\bar{\boldsymbol{x}}^{s,2}_{k+1}]^{T}\bar{\boldsymbol{x}}^{s,2}_{k+1}\nonumber\\
	&&\quad
	+[\bar{\boldsymbol{x}}^{s+1,2}_{k+1}-\bar{\boldsymbol{x}}^{s,2}_{k+1}]^{T}[\bar{\boldsymbol{x}}^{s+1,2}_{k+1}-\bar{\boldsymbol{x}}^{s,2}_{k+1}]. \label{V1}
\end{eqnarray}
Now we give upper bounds of each term of the right side of~\eqref{V1} based on~\eqref{311}. The first term is written as
\begin{equation}
[\bar{\boldsymbol{x}}^{s+1,1}_{k+1}-\bar{\boldsymbol{x}}^{s,1}_{k+1}]^{T}[\bar{\boldsymbol{x}}^{s+1,1}_{k+1}-\bar{\boldsymbol{x}}^{s,1}_{k+1}]\leq\|\left(M_1^{T} \otimes I_N\right) R\Theta^{s}_{k+1}\|^2. \nonumber
\end{equation}
Observing $\lambda_{2}\leq \|M_{2}^{T}LM_{2}\| \leq \lambda_{N}$, we have
\begin{eqnarray}
	&&\quad [\bar{\boldsymbol{x}}^{s+1,2}_{k+1}-\bar{\boldsymbol{x}}^{s,2}_{k+1}]^{T}[\bar{\boldsymbol{x}}^{s+1,2}_{k+1}-\bar{\boldsymbol{x}}^{s,2}_{k+1}] \nonumber\\
	&&\leq2\|\alpha(M_2^{T} L M_2\otimes I_N) \bar{\boldsymbol{x}}^{s,2}_{k+1}\|^2\!+\!2\|\left(M_2^{T} \otimes I_N\right) R  \Theta^{s}_{k+1}\|^2,\nonumber\\
	&&\leq 2\alpha^2\lambda_{N}^2\|\bar{\boldsymbol{x}}^{s,2}_{k+1}\|^2\!+\!2\|\left(M_2^{T} \otimes I_N\right) R \Theta^{s}_{k+1}\|^2.\nonumber
\end{eqnarray}
Applying the same argument, we have
\begin{eqnarray}
		&&2[\bar{\boldsymbol{x}}^{s+1,1}_{k+1}\!-\!\bar{\boldsymbol{x}}^{s,1}_{k+1}]^{T}\bar{\boldsymbol{x}}^{s,1}_{k+1}\!=\!-2(\bar{\boldsymbol{x}}^{s,1}_{k+1})^{T}\left(M_1^{T} \otimes I_N\right) R\Theta^{s}_{k+1},\nonumber\\
		&&2[\bar{\boldsymbol{x}}^{s+1,2}_{k+1}-\bar{\boldsymbol{x}}^{s,2}_{k+1}]^{T}\bar{\boldsymbol{x}}^{s,2}_{k+1}\leq -2\alpha \lambda_{2}\|\bar{\boldsymbol{x}}^{s,2}_{k+1}\|^2-2(\bar{\boldsymbol{x}}^{s,2}_{k+1})^{T}\nonumber\\
		&& \quad \quad \quad \quad \quad \quad \quad \quad \quad \quad \quad \quad \quad\left(M_2^{T} \otimes I_N\right)R \Theta^{s}_{k+1}.\label{3333}
\end{eqnarray}
By inserting~\eqref{3333}  into~\eqref{V1}, applying the relationships
\begin{equation}
\begin{cases}
\tilde{\boldsymbol{x}}^{s}_{k+1}=(\bar{\boldsymbol{x}}^{s,1}_{k+1})^{T}[M_1^{T} \otimes I_N]+(\bar{\boldsymbol{x}}^{s,2}_{k+1})^{T}[M_2^{T} \otimes I_N],\\
\|R \Theta^{s}_{k+1}\|^2=\|\left(M_1^{T} \otimes I_N\right)R \Theta^{s}_{k+1}\|^2\\
\quad \quad\quad\quad\quad\quad\quad\quad+\|\left(M_2^{T} \otimes I_N\right) R \Theta^{s}_{k+1}\|^2,\nonumber
\end{cases}
\end{equation}
and $\|M_{2}\|=1$, it derives that
\begin{eqnarray}
&&V^{s+1}_{k+1}-V^{s}_{k+1}\leq -2(\tilde{\boldsymbol{x}}^{s}_{k+1})^{T}R \Theta^{s}_{k+1}+2\|R \Theta^{s}_{k+1}\|^2\nonumber\\
&&\quad \quad \quad \quad \quad \quad \quad  +2\alpha^2\lambda_{N}^2\|\bar{\boldsymbol{x}}^{s,2}_{k+1}\|^2-2\alpha \lambda_{2}\|\bar{\boldsymbol{x}}^{s,2}_{k+1}\|^2.\label{dotV}
\end{eqnarray}
We split $\tilde{\boldsymbol{x}}^{s}_{k+1}$ into $\tilde{\boldsymbol{x}}^{s,1}_{k+1}+\tilde{\boldsymbol{x}}^{s,2}_{k+1}$ to estimate
\begin{eqnarray}
&&\quad -(\tilde{\boldsymbol{x}}^{s}_{k+1})^{T}R \Theta^{s}_{k+1}\nonumber\\
&&=-(\tilde{\boldsymbol{x}}^{s,1}_{k+1})^{T}R[\boldsymbol{F}^{t_{k+1}}(\tilde{\boldsymbol{x}}^{s,1}_{k+1}\!\!+\!\boldsymbol{x}^*(t_{k+1}))\!\!-\!\boldsymbol{F}^{t_{k+1}}\left(\boldsymbol{x}^*(t_{k+1})\right)]\nonumber\\
&&\quad-(\tilde{\boldsymbol{x}}^{s,2}_{k+1})^{T}R[\boldsymbol{F}^{t_{k+1}}(\tilde{\boldsymbol{x}}^{s,1}_{k+1}\!\!+\!\boldsymbol{x}^*(t_{k+1}))\!\!-\!\boldsymbol{F}^{t_{k+1}}\left(\boldsymbol{x}^*(t_{k+1})\right)]\nonumber\\
&&\quad-(\tilde{\boldsymbol{x}}^{s,1}_{k+1})^{T}R[\boldsymbol{F}^{t_{k+1}}(\tilde{\boldsymbol{x}}^{s,1}_{k+1}+\tilde{\boldsymbol{x}}^{s,2}_{k+1}+\boldsymbol{x}^*(t_{k+1}))\nonumber\\
&&\quad \quad \quad \quad\quad \quad\quad \quad\quad \quad\quad  -\boldsymbol{F}^{t_{k+1}}(\tilde{\boldsymbol{x}}^{s,1}_{k+1}+\boldsymbol{x}^*(t_{k+1}))]\nonumber\\
&&\quad-(\tilde{\boldsymbol{x}}^{s,2}_{k+1})^{T}R[\boldsymbol{F}^{t_{k+1}}(\tilde{\boldsymbol{x}}^{s,1}_{k+1}+\tilde{\boldsymbol{x}}^{s,2}_{k+1}+\boldsymbol{x}^*(t_{k+1}))\nonumber\\
&&\quad \quad \quad \quad\quad \quad\quad \quad\quad \quad\quad  -\boldsymbol{F}^{t_{k+1}}(\tilde{\boldsymbol{x}}^{s,2}_{k+1}+\boldsymbol{x}^*(t_{k+1}))]
.\label{dotV1}
\end{eqnarray}
As for $\boldsymbol{F}^{t_{k+1}}(1_{N}\otimes y)=F^{t_{k+1}}(y)$ for any $y\in{\mathbb{R}^{N}}$, it follows by the strong monotonically of $F^{t_{k+1}}$ that
\vspace{-0.5em}
\begin{eqnarray}
&&(\tilde{\boldsymbol{x}}^{s,1}_{k+1})^{T}R[\boldsymbol{F}^{t_{k+1}}(\tilde{\boldsymbol{x}}^{s,1}_{k+1}+\boldsymbol{x}^*(t_{k+1}))-\boldsymbol{F}^{t_{k+1}}\left(\boldsymbol{x}^*(t_{k+1})\right)]\nonumber\\
&&=\frac{\varepsilon(\bar{\boldsymbol{x}}^{s,1}_{k+1})^{T}}{\sqrt{N}}[F^{t_{k+1}}(\frac{(\bar{\boldsymbol{x}}^{s,1}_{k+1})^{T}}{\sqrt{N}}+x^*(t_{k+1}))\nonumber\\
&&\quad \quad \quad \quad \quad\quad \quad -F^{t_{k+1}}( x^*(t_{k+1}))]\geq \frac{\varepsilon\mu}{N}\|\bar{\boldsymbol{x}}^{s,1}_{k+1}\|^2,\label{dotV11}
\end{eqnarray}
where $(1_{N}^{T}\otimes I_{N})R=I_{N}$ and $(\tilde{\boldsymbol{x}}^{s,1}_{k+1})^{T}R=\frac{\varepsilon(\bar{\boldsymbol{x}}^{s,1}_{k+1})^{T}}{\sqrt{N}}$ are used. Since $\|R^{T}\tilde{\boldsymbol{x}}^{s,2}_{k+1}\|=\varepsilon\|\tilde{\boldsymbol{x}}^{s,2}_{k+1}\|=\varepsilon\|\bar{\boldsymbol{x}}^{s,2}_{k+1}\|$, 
\begin{eqnarray}
&&(\tilde{\boldsymbol{x}}^{s,2}_{k+1})^{T}R[\boldsymbol{F}^{t_{k+1}}(\tilde{\boldsymbol{x}}^{s,1}_{k+1}+\boldsymbol{x}^*(t_{k+1}))-\boldsymbol{F}^{t_{k+1}}\left(\boldsymbol{x}^*(t_{k+1})\right)]\nonumber\\
&&\leq \frac{\varepsilon\theta^{t_{k+1}}}{\sqrt{N}}\|\bar{\boldsymbol{x}}^{s,1}_{k+1}\|\|\bar{\boldsymbol{x}}^{s,2}_{k+1}\|.\label{dotV12}
\end{eqnarray}
Similar operations are used to obtain that
\vspace{-0.5em}
\begin{flalign}
-[\tilde{\boldsymbol{x}}^{s}_{k+1}]^{T}R \Theta^{s}_{k+1} &\leq \frac{2\varepsilon\theta^{t_{k+1}}}{\sqrt{N}}\|\bar{\boldsymbol{x}}^{s,1}_{k+1}\|\|\bar{\boldsymbol{x}}^{s,2}_{k+1}\|\nonumber\\
& \quad \quad  +\varepsilon\theta^{t_{k+1}}\|\bar{\boldsymbol{x}}^{s,2}_{k+1}\|^2\!-\!\frac{\varepsilon\mu^{t_{k+1}}}{N}\|\bar{\boldsymbol{x}}^{s,1}_{k+1}\|^2.\nonumber
\end{flalign}
Then the difference of $V_{k+1}$ in~\eqref{dotV} is reorganized as
\vspace{-0.5em}
\begin{eqnarray}
&&\quad V^{s+1}_{k+1}-V^{s}_{k+1}\nonumber\\
&&\leq \frac{4\varepsilon\theta^{t_{k+1}}}{\sqrt{N}}\|\bar{\boldsymbol{x}}^{s,1}_{k+1}\|\|\bar{\boldsymbol{x}}^{s,2}_{k+1}\|\!-\!\frac{2\varepsilon\mu^{t_{k+1}}}{N}\|\bar{\boldsymbol{x}}^{s,1}_{k+1}\|^2\!+\!2\|R \Theta^{s}_{k+1}\|^2\nonumber\\
&&\quad \quad \quad \quad\quad \quad +(2\varepsilon\theta^{t_{k+1}}+2\alpha^2\lambda_{N}^2-2\alpha \lambda_{2})\|\bar{\boldsymbol{x}}^{s,2}_{k+1}\|^2.\label{dotVV}
\end{eqnarray}
We further bound $\|R \Theta^{s}_{k+1}\|$ by
\vspace{-0.5em}
\begin{eqnarray}
&&R \Theta^{s}_{k+1}=R[\boldsymbol{F}^{t_{k+1}}(\tilde{\boldsymbol{x}}^{s,1}_{k+1}+\tilde{\boldsymbol{x}}^{s,2}_{k+1}+\boldsymbol{x}^*(t_{k+1}))\nonumber\\
&&~-\boldsymbol{F}^{t_{k+1}}(\tilde{\boldsymbol{x}}^{s,1}_{k+1}\!+\!\boldsymbol{x}^*(t_{k+1}))]\!+\! R[\boldsymbol{F}^{t_{k+1}}(\tilde{\boldsymbol{x}}^{s,1}_{k+1}\!+\!\boldsymbol{x}^*(t_{k+1}))\nonumber\\
&&~-\boldsymbol{F}^{t_{k+1}}\left(\boldsymbol{x}^*\right)]\leq \varepsilon \theta^{t_{k+1}} \|\bar{\boldsymbol{x}}^{s,2}_{k+1}\|+\varepsilon \theta^{t_{k+1}} \|\bar{\boldsymbol{x}}^{s,1}_{k+1}\|.
\end{eqnarray}
To sum up, we conclude that
\vspace{-0.5em}
\begin{flalign}
V^{s+1}_{k+1}\!-\!V^{s}_{k+1}&\leq\! -2[\|\bar{\boldsymbol{x}}^{s,1}_{k+1}\|~~\|\bar{\boldsymbol{x}}^{s,2}_{k+1}\|]A_{\alpha,\varepsilon}\left[\begin{matrix}
		\|\bar{\boldsymbol{x}}^{s,1}_{k+1}\| \\
		\|\bar{\boldsymbol{x}}^{s,2}_{k+1}\|
		\end{matrix}\right],\nonumber\\
	&\leq -\nu V^{s}_{k+1},
\end{flalign}
which implies that
\vspace{-0.5em}
\begin{eqnarray}
	\|\hat{\boldsymbol{x}}^{s+1}_{k+1}-\boldsymbol{x}^{*}(t_{k+1})\|^2 \leq (1-\nu) \|\hat{\boldsymbol{x}}^{s}_{k+1}-\boldsymbol{x}^{*}(t_{k+1})\|^2.\label{s0}
\end{eqnarray}
Since $\hat{\boldsymbol{x}}^{s}_{k+1}$ is initialized by the predicted variable $\boldsymbol{x}_{k+1|k}$ and the corrected variable $\boldsymbol{x}_{k+1}=\hat{\boldsymbol{x}}^{\tau}_{k+1}$, considering the relation in~\eqref{s0} between two consecutive iterates of the sequence $\hat{\boldsymbol{x}}^{s}_{k+1}$, we obtain~\eqref{correctionerror}.
\vspace{-0.5em}
\bibliographystyle{apacite}    
\bibliography{mybibfile_nonlinear}           


\end{document}